\documentclass[AMA,STIX1COL]{WileyNJD-v2}

\articletype{Article Type}%
\received{ }
\revised{ }
\accepted{ }

\begin{document}

\title{Lie symmetries and travelling wave solutions of the nonlinear waves in the inhomogeneous Fisher-Kolmogorov equation.}

\author[1]{M.S. Bruz\'on*}

\author[1]{T.M. Garrido}

\author[1]{E. Recio}

\author[1]{R. de la Rosa}

\authormark{M.S. Bruz\'on \textsc{et al}}

\address{\orgdiv{Department of Mathematics}, \orgname{University of C\'{a}diz}. \orgaddress{PO.BOX 40, 11510 Puerto Real, \country{Spain}}}

\corres{*M.S. Bruz\'on. \email{m.bruzon@uca.es}}


\abstract[Summary]{In this work we consider a Fisher-Kolmogorov equation depending on two exponential functions of the spatial variables. We study this equation from the point of view of symmetry reductions in partial differential equations. Through two-dimensional abelian subalgebras, the equation is reduced to ordinary differential equations. New solutions have been derived and interpreted.}


\jnlcitation{\cname{%
\author{M.S. Bruz\'on},
\author{T.M. Garrido},
\author{E. Recio}, and
\author{R. de la Rosa}} (\cyear{2019}),
\ctitle{Lie symmetries of the nonlinear waves in the inhomogeneous Fisher-Kolmogorov equation}, \cjournal{Math Meth Appl Sci.}, \cvol{2019;00:1--9}.}

\maketitle

\section{Introduction}
In 1930, the Fisher-Kolmogorov (FK) equation was proposed for populations dynamics. It shows spread of an advantageous gene in a population. The existence of solutions and travelling waves was demonstrated in 1937. The analysis and study of the Fisher equation is used to model heat and reaction-diffusion problems applied to mathematical biology, physics, astrophysics, chemistry, genetics, bacterial growth problems as well as development and growth of solid tumours. There are several variations of the FK equation, let us show few of them.\\

One of the well-known models is the Fisher-Kolmogorov equation or the Fisher-Kolmogorov-Petrovskii-Piskunov equation, which was proposed, on the one hand, by Fisher \cite{FK1} who described the spread of an advantageous gene in a population; and on the other hand, by Kolmogorov et al. \cite{FK2} who additionally obtained the basic analytical results for this equation. This two-dimensional FK equation has the following form:
\begin{equation}\label{Fis}
u_t= au(1-u)+b\,\Delta u,
\end{equation}
where $u$ is the frequency of the mutant gene, $a$ is an intrinsic coefficient, $\Delta = \partial^2/\partial x^2+\partial^2/\partial y^2$ is the Laplace operator, and $b$ is the diffusion constant.\\

\noindent In \cite{GBR}, the following generalized Fisher equation was considered
\begin{equation}\label{eqF1}
\begin{array}{l}
u_{t}=f(u)+(g(u)u_x)_x,
\end{array}
\end{equation}
where the diffusivity depends on an arbitrary function $g(u)$, being $x$ and $t$ the independent variables, and $f(u)$ an arbitrary function. For equation (\ref{eqF1}), some non-trivial conservation laws have been obtained.\\

\noindent Another generalized Fisher equation was considered in \cite{RBG}
\begin{equation}\label{eqF2}
\begin{array}{l}
u_{t}=f(u)+\frac{1}{x}(xg(u)u_x)_x,
\end{array}
\end{equation}
where $x$ is used as the radial variable, $g(u)$ is the diffusion coefficient, and $f(u)$ is an arbitrary function. For equation (\ref{eqF2}) the Lie group classification was obtained, and all the reductions were derived from the optimal system of subalgebras. Moreover, as some of the reduced equations admit Lie symmetries, further reductions were derived. In addition, some exact wave solutions were obtained by using a direct method.\\

\noindent The so called inhomogeneous FK equation,
\begin{equation}\label{eq11}
u_t=\left[a(x)-c(x)u\right]u^n+\left[b(x)u_x\right]_x,
\end{equation}
was considered in \cite{Shi}. Equation (\ref{eq11}) models the population grow in an heterogeneous environment. Here, the inhomogeneity appears in the diffusion coefficient $b(x)$ and the intrinsic growth rate $a(x)$. The simplest function to consider as inhomogeneity is a periodic function. As a result, the spread of population waves located in a bounded region leads to formation and development of periodic patterns moving with a constant velocity. It is observed that the restriction by a finite continuous function indicates that it is quite natural to choose the simplest form for the function $c(x)=\frac{1}{\cos(x)}$.\\

Additionally, in \cite{main} was studied how inhomogeneities in equation (\ref{eq11}) influence on the population invasions and the velocity of periodic waves. This case was considered with $a$, $b$ constants and $n=1$. In order to study the previous systems, they considered the equation in the general form
\begin{equation}\label{eq11K}
u_t=au-\frac{u^2}{K(x)}+b u_{xx}.
\end{equation}
The solutions for (\ref{eq11K}) obtained in \cite{Shi,main} admit negative values for the function $u$. However, in the context of the population dynamics this function should be non-negative, i.e. in the simplest case it should have the form $u(x, t) \propto 1+\cos(x)$ instead of $\cos (x)$. Therefore, to find the exact solution for such ``real'' waves they considered a more general variant of the FK equation with inhomogeneous reaction and diffusion terms. So, equation (\ref{eq11}) was considered  with $n=1$. Hence, all inhomogeneities have similar non-negative dependencies. Particular cases were studied to obtain exact analytical solutions. For example, for $a(x)=3b_0(1+\cos(x))$, $c(x)=\frac{c_0}{1+\cos(x)}$, $b(x)=b_0(1+\cos(x))$ and $n=1$, with $a_0$, $b_0$, $c_0$ arbitrary constants, the authors showed that equation (\ref{eq11}) admits the solution
\begin{equation}\label{solf}
u(x,t)=\frac{3b_0(1+\cos(x))}{e^{k-3b_0t}+c_0},
\end{equation} where $k$ is an arbitrary constant.\\

\noindent In this paper, we analyse a two-dimensional generalization of equation (\ref{eq11K}) given by
\begin{equation}\label{edp}
u_t=au - M(x)N(y)\ u^2 + b\left(u_{xx} + u_{yy}\right),
\end{equation}
where $b$ is the diffusion coefficient, $a$ is the intrinsic growth rate, and the inhomogeneity can be factorized in the form $M(x)N(y)$ \cite{main}. For computational aspects, $M(x), N(y)$ have been considered in the numerator.
The variable $u$ represents the population density, so $u$ is non-negative. The term $au$ represents a source if $a>0$ or a sink if $a<0$, and $M(x)N(y)$ will be considered non-negative since $-M(x)N(y)u^2$ represents a sink.
Equation (\ref{edp}) is called the generalized (2+1)-dimensional inhomogeneous Fisher-Kolmogorov equation. Some special solutions of (\ref{edp}) have been obtained in \cite{main}. \\

Lie symmetry method allows us to reduce the number of independent variables of a given partial differential equation (PDE). 
In particular, it allows to reduce a (2+1)-dimensional PDE into a (1+1)-dimensional PDE.
When a PDE in $2+1$ dimensions admits a 2-dimensional subalgebra of symmetries,
then it is possible to reduce the PDE to an ordinary differential equation (ODE). 
Therefore, the Lie method is an important and efficient tool for analysing differential equations \cite{bluman,olver}. 
Applications of the Lie method include determining invariant solutions, constructing maps between equivalent equations, determining conserved quantities when the PDE admits a Lagrangian formulation, and so on.
Symmetries have been obtained for two dimensional models to look for exact solutions in \cite{recio,khaliq,kumar} with effective results.\\

The main purpose of the present paper is to obtain new solutions for equation (\ref{edp}). 
We focus our attention on the determination of invariant solutions, i.e. solutions which are invariant under the action of Lie symmetry groups. 
First, we carry out a classification of Lie point symmetries admitted by equation (\ref{edp}). 
For simplicity, it is assumed that $M(x)$ are $N(y)$ are exponential functions. 
A complete study of the symmetries for arbitrary $M(x)$ are $N(y)$ will appear elsewhere.
In addition, we present the symmetry transformation groups corresponding to the admitted point symmetries.
Next, we find some exact solutions for equation (\ref{edp}), taking into account two-dimensional subalgebras admitted by the equation, and we interpret those solutions. Finally, we give some conclusions.

\section{Lie symmetries}\label{sec:symm}

We apply the Lie group method of infinitesimal transformations to equation (\ref{edp}) in the case $M(x)=m_1 e^{px}$, $N(y)=n_1 e^{qy}$. Following references \cite{bluman,olver}, we consider a one-parameter Lie group of infinitesimals transformations acting on independent and dependent variables given by
$$\begin{array}{rcl}
\displaystyle \tilde t =& t &+\ \varepsilon \tau (t,x,y,u) + \mathcal{O} (\varepsilon^2), \\ \vspace{0.1cm}
\displaystyle \tilde x =& x &+\ \varepsilon \xi^x (t,x,y,u) + \mathcal{O} (\varepsilon^2),\\ \vspace{0.1cm}
\displaystyle \tilde y =& y &+\ \varepsilon \xi^y (t,x,y,u) + \mathcal{O} (\varepsilon^2), \\ \vspace{0.1cm}
\displaystyle \tilde u =& u &+\ \varepsilon \phi (t,x,y,u) + \mathcal{O} (\varepsilon^2),
\end{array}$$
where $\varepsilon$ is the group parameter and the associated vector field takes the following form
\begin{equation}\label{vect}
X=\tau (t,x,y,u) \partial_t  + \xi^x (t,x,y,u)\partial_x+\xi^y (t,x,y,u)\partial_y+  \phi (t,x,y,u) \partial_u.
\end{equation}

\noindent The invariance criterion implies
\begin{equation} \label{edet}
\mbox{pr}^{(2)}X \Delta = 0, \quad \mbox{where} \quad \Delta=0,
\end{equation}
with $\Delta=u_t-au+ M(x) N(y) u^2 -b\left(u_{xx}+u_{yy}\right)$, and where $\mbox{pr}^{(2)}X$ is the second prolongation of the vector field (\ref{vect}) defined by
$$\mbox{pr}^{(2)}X= X +\phi_t \frac{\partial}{\partial u_t} + \phi_x \frac{\partial}{\partial u_x} + \phi_y \frac{\partial}{\partial u_y} + \phi_{xx} \frac{\partial}{\partial u_{xx}}+ \phi_{xy} \frac{\partial}{\partial u_{xy}}+ \phi_{yy} \frac{\partial}{\partial u_{yy}},$$
with coefficients
$$ \begin{array}{ll}
\phi_t =& D_t(\eta)- u_tD_t(\tau) - u_xD_t(\xi^x) - u_y D_t(\xi^y),\\
\phi_x =& D_x(\eta)- u_tD_x(\tau) - u_xD_x(\xi^x) - u_y D_x(\xi^y),\\
\phi_y =& D_y(\eta)- u_tD_y(\tau) - u_xD_y(\xi^x) - u_y D_y(\xi^y),\\
\phi_{xx} =& D_x(\phi_x)- u_{tx}D_x(\tau) - u_{xx}D_x(\xi^x) - u_{yx} D_x(\xi^y),\\
\phi_{xy} =& D_y(\phi_x)- u_{ty}D_y(\tau) - u_{xy}D_y(\xi^x) - u_{yy} D_y(\xi^y),\\
\phi_{yy} =& D_y(\phi_y)- u_{ty}D_y(\tau) - u_{xy}D_y(\xi^x) - u_{yy} D_y(\xi^y),
\end{array}
$$
where $D_t, D_x, D_y$ are the total derivatives of $t$, $x$ and $y$, respectively.\\ 

Equation (\ref{edet}) splits with respect to the derivatives of $u$. This yields an overdetermined linear system of equations for the infinitesimals $\tau (t,x,y,u)$, $\xi^x (t,x,y,u)$, $\xi^y (t,x,y,u)$, $\phi (t,x,y,u)$ along with functions $M(x), N(y)$ and the parameters $a$ and $b$, with the conditions $M(x)\neq0$, $N(y)\neq0$, $a\neq0$ and $b\neq0$. Equation (\ref{edp}) has several obvious discrete equivalence transformations. In particular, it should be noted that it admits the discrete equivalence transformation
\begin{equation}\label{disctr}
t \longrightarrow \tilde{t}, \qquad x \longrightarrow \tilde{y}, \qquad y \longrightarrow \tilde{x}, \qquad u \longrightarrow \tilde{u},  \qquad M \longrightarrow \tilde{N}, \qquad N \longrightarrow \tilde{M},
\end{equation}
i.e., a transformation under which $x$ and $y$ are interchanged along with $M(x)$ and $N(y)$.\\

From the overdetermined linear system of equations we deduce that $\tau=\tau (t)$, $\xi^x =\xi^x(t,x,y)$, $\xi^y=\xi^y (t,x,y)$, $\phi= \phi (t,x,y,u)$, where $\tau, \xi^x, \xi^y$, $\phi$, $M$ and $N$, must satisfy the system
\begin{align}
\phi_{uu}=0,\\
\tau_t = 2 \xi^x_x,\\
\xi^y_x=-\xi^x_y,\\
\xi^y_y = \xi^x_x,\\
\phi_{xu} = \frac{1}{2 b} \left( -\xi^x_t + b \xi^x_{xx} + b \xi^x_{yy} \right), \\
\phi_{yu} = \frac{-1}{2b} \xi^y_t,\\
\phi_t - b \left( \phi_{xx} + \phi_{yy} \right) + \left( a -MN u \right) u \phi_u + \left( -a+2 MN u \right) \phi + \left( -2 au + 2MN u^2 \right) \xi^x_x + M_x N u^2 \xi^x + MN_yu^2 \xi^y=0.
\end{align}

We solve this system for the case $M(x)=m_1 e^{px}$ and $N(y)=n_1 e^{qy}$ by using Maple "rifsimp" and "pdsolve" commands. Thus, we obtain the following result.

\begin{theorem}
The complete classification of the point symmetries admitted by the (2+1)-dimensional inhomogeneous Fisher-Kolmogorov equation (\ref{edp}) with $M(x)=m_1 e^{px}$, $N(y)=n_1 e^{qy}$ is given by:
\begin{itemize}
\item Case 1. For $a$, $b$, $m_1$, $n_1$ non-zero arbitrary parameters and $p$, $q$ arbitrary parameters, the infinitesimal generators are
\begin{align}
X_1= \partial_t, \label{symm1}\\
X_2= \partial_x - p u \partial_u, \label{symm2} \\
X_{2'}= \partial_y -q u \partial_u, \label{symm2'} \\
X_3= \left( -y + 2b q t \right) \partial_x + \left( x- 2b pt \right) \partial_y + \left( py-qx \right) u \partial_u.
\label{symm3}
\end{align}

\item Case 2. For $a=-b \left( p^2 + q^2 \right)$, the infinitesimal generators are $X_1, X_2, X_{2'}, X_3$ and
\begin{equation}\label{symm4}
X_4=pqt \partial_t + \left( \frac{1}{2}pq x + b p^2 q t - \frac{1}{2}q \right) \partial_x + \left( \frac{1}{2}pq y + b pq^2t - \frac{1}{2}p \right) \partial_y - \left( b \left(p^2 + q^2 \right)pq t + \frac{1}{2} p^2qx + \frac{1}{2}pq^2y \right) u \partial_u.
\end{equation}

\end{itemize}

\end{theorem}

\noindent Each symmetry (\ref{vect}) generates a transformation obtained by solving the following system of ODEs
$$\begin{array}{rcl}
\dfrac{\partial \tilde{t}}{\partial \epsilon}=\tau(\tilde{t},\tilde{x},\tilde{y},\tilde{u}),\qquad
\dfrac{\partial \tilde{x}}{\partial \epsilon}=\xi^x(\tilde{t},\tilde{x},\tilde{y},\tilde{u}),\qquad
\dfrac{\partial \tilde{y}}{\partial \epsilon}=\xi^y(\tilde{t},\tilde{x},\tilde{y},\tilde{u}),\qquad
\dfrac{\partial \tilde{u}}{\partial \epsilon}=\phi(\tilde{t},\tilde{x},\tilde{y},\tilde{u}),
\end{array}$$
satisfying the initial conditions
$$\tilde{t}|_{\epsilon=0}=t, \qquad \tilde{x}|_{\epsilon=0}=x, \qquad \tilde{y}|_{\epsilon=0}=y, \qquad \tilde{u}|_{\epsilon=0}=u,$$
with $\epsilon$ the group parameter. Point symmetries (\ref{symm1})-(\ref{symm4}) produce the following one-parameter symmetry transformation groups:

$$\begin{array}{l}
(\tilde{t},\tilde{x},\tilde{y},\tilde{u})_1=
(t+\epsilon,x,y,u), 
\qquad 
\mbox{time translation,}
\\ \\
(\tilde{t},\tilde{x},\tilde{y},\tilde{u})_2=
(t,x+\epsilon,y,e^{-p \epsilon} u), 
\qquad 
\mbox{scaling in } u \mbox{ combined with translation in } x,
\\ \\
(\tilde{t},\tilde{x},\tilde{y},\tilde{u})_{2'}=
(t,x,y+\epsilon,e^{-q \epsilon} u), 
\qquad 
\mbox{scaling in } u \mbox{ combined with translation in } y ,
\\ \\
(\tilde{t},\tilde{x},\tilde{y},\tilde{u})_{3}=
\Big(t,
(x-2b pt)\cos(\epsilon)+(2b qt-y)\sin(\epsilon)+2b pt,
(y-2b qt)\cos(\epsilon)+(x-2b pt)\sin(\epsilon)+2b qt,
\\
\qquad\qquad\qquad\quad
\exp\big((py-qx)\sin(\epsilon)+(2b p^2 t-px+b q^2t-qy)(\cos(\epsilon)-1)\big)u
\Big),
\\
\qquad\qquad\qquad\qquad\qquad\qquad\qquad\quad
\mbox{dilation combined with a rotation in the plane } (x,y),
\\ \\
(\tilde{t},\tilde{x},\tilde{y},\tilde{u})_{4}=
\Big(
e^{\epsilon}t,
(e^{\epsilon/2}-1)(2b pe^{\epsilon/2} t -\frac{1}{p})+e^{\epsilon/2} x,
(e^{\epsilon/2}-1)(2b qe^{\epsilon/2} t -\frac{1}{q})+e^{\epsilon/2}y,
\\
\qquad\qquad\qquad\quad
\exp\big((2b(p^2+q^2)t-px-qy+2)e^{\epsilon/2}
-2b(p^2+q^2)e^{\epsilon}t+px+qy-\epsilon-2\big)u
\Big),
\\
\qquad\qquad\qquad\qquad\qquad\qquad\qquad\quad
\mbox{dilation combined with a scaling}.
\end{array}
$$

\section{Reductions and some exact solutions}

Point symmetries of equation (\ref{edp}) can be used to transform it into equations with fewer number of independent variables. In order to determine exact solutions of equation (\ref{edp}), we focus our attention in two-dimensional abelian subalgebras. Thus, equation (\ref{edp}) can be reduced into a second-order nonlinear ODE.


\subsection{Reduction under the two-dimensional abelian subalgebra $\left\{ X_1-\alpha X_2, X_{2'}+\beta X_2 \right\}$}

Consider $M(x)=m_1 e^{p x}$ and $N(y)=n_1 e^{q y}$ and the two-dimensional abelian subalgebra
\begin{equation}\label{algebra} 
\mathscr{A}: \mbox{span} \left( X_1-\alpha X_2, X_{2'}+\beta X_2 \right).
\end{equation}

By taking into account the two-dimensional abelian subalgebra (\ref{algebra}) one can reduce (\ref{edp}) into a nonlinear second-order ODE. Let $X=X_1-\alpha X_2$ and $Y=X_{2'}+\beta X_2$. Consequently, one looks for two independent invariants $z$ and $w$ satisfying
\begin{equation}
\label{conditions}
\begin{array}{l}
X\, z=0, \qquad X \, w=0, \\
Y\, z=0, \qquad Y\, w=0.
\end{array}
\end{equation}

\noindent From (\ref{conditions}) we obtain
$$ \displaystyle z=-x+\beta y-\alpha t, \qquad w= e^{px+qy} u.$$

\noindent Thus, 
\begin{equation}\label{grinvsol}
u= e^{-(p x+q y)} w \left(-x+\beta y-\alpha t\right),
\end{equation}

\noindent is a group-invariant solution. The group-invariant solution (\ref{grinvsol}) reduces equation (\ref{edp}) into a second-order ODE, given by
\begin{equation}
\label{ode1} 
w''= \frac{2b (\beta q-p)-\alpha}{ b (\beta^2+1) } w'+ \frac{m_1 n_1}{b (\beta^2+1) } w^2-\frac{b ( p^2+q^2)+a}{ b (\beta^2+1) } w.
\end{equation}

\noindent In \cite {NK}, the author obtained the general solution of a second-order ODE of the form
\begin{equation}\label{eqK}  w''=\dfrac{1}{\lambda} \left( \mu w'+\dfrac{1}{2 } w^2 -\omega w -c_0 \right), \end{equation} with $c_0$ an arbitrary constant and $\lambda$, $\mu$ and $\omega$ satisfying the constraint
\begin{equation}
\label{constraint} \omega= \frac{6 \mu^2}{25 \lambda}.
\end{equation}
Equation (\ref{ode1}) can be written in the same form that equation (\ref{eqK}) by considering 
\begin{eqnarray}
\label{var1}\lambda&=&\frac{b(\beta^2+1)}{2 m_1 n_1},\\
\label{var2}\mu&=& \frac{2b (\beta q-p)-\alpha}{2 m_1 n_1},\\
\label{var3}\omega&=& \frac{b ( p^2+q^2)-a}{2 m_1 n_1},\\
\label{var4}c_0&=& 0,
\end{eqnarray} 
where the constraint (\ref{constraint}) takes the form
\begin{equation}
\label{constraintode} 
25 b ( b ( p^2+q^2)+a) (\beta^2+1)= 6 \left( 2 b(p-\beta q)+\alpha \right)^2.
\end{equation}

\noindent So, the general solution of equation (\ref{ode1}) is given by
\begin{equation}\label{sode1} 
w(z)=  \frac{3 b (\beta^2+1)}{2 m_1 n_1}\exp\left(\frac{2 (2 b(\beta q-p)-\alpha)z}{5 b(\beta^2+1)}\right) \mathcal{P}\left(\frac{5 b(\beta^2+1)}{4 b(\beta q-p)-2 \alpha}\exp\left(\frac{ (2 b(\beta q-p)-\alpha)z}{5 b(\beta^2+1)}\right)+c_1, 0, c_2\right),
\end{equation} 
in terms of the Weierstrass elliptic function with invariants $g_2=0$ and $g_3=c_2$, where $c_1, c_2$ are arbitrary constants.
\\

Finally, by undoing the change of variables (\ref{grinvsol}), we obtain the following solution of the (2+1)-dimensional inhomogeneous Fisher-Kolmogorov equation (\ref{edp}):
\begin{equation}\label{solA}
\begin{array}{l}
\displaystyle u(t,x,y)= 
\frac{3 b (\beta^2+1)}{2 m_1 n_1}
\exp\big(-(p x+q y)\big)
\exp\left(\frac{2 (2 b(\beta q-p)-\alpha)(\beta y-x-\alpha t)}{5 b(\beta^2+1)}\right) 
\\
\qquad\qquad\qquad\quad
\displaystyle \mathcal{P}\left(\frac{5 b(\beta^2+1)}{4 b(\beta q-p)-2 \alpha}\exp\left(\frac{ (2 b(\beta q-p)-\alpha)(\beta y-x-\alpha t)}{5 b(\beta^2+1)}\right)+c_1, 0, c_2\right).
\end{array}
\end{equation}
\\

\noindent 
Solution (\ref{solA}) can be interpreted
by considering the geometrical form of the travelling wave variable 
\begin{equation}
z=-x+\beta y - \alpha t = |{\mathbf k}|( \hat{\mathbf k}\cdot (x,y) - \alpha t ),
\end{equation}
where ${\mathbf k}=(-1,\beta)$ is a constant vector in the $(x,y)$-plane, $|{\mathbf k}|$ is its modulus, and 
where $\hat{\mathbf k}$ is the unit vector
\begin{equation}
\hat{\mathbf k} = (\cos\phi,\sin\phi),
\quad
\tan\phi = -\beta,
\end{equation}
giving the direction of propagation of the solution,
with speed
\begin{equation}
c = \alpha/|{\mathbf k}|, 
\quad
|{\mathbf k}|^2 = 1+\beta^2.
\end{equation}
Since the coefficient of $u^2$ in equation (\ref{edp}) is given by $M(x)N(y)=m_1 n_1 e^{px+qy}$,
to interpret the solution,
it is natural to consider $px+qy= (p,q) \cdot (x,y) = {\mathbf l} \cdot {\mathbf x} = |{\mathbf l}| \hat{\mathbf l} \cdot {\mathbf x}$, where $\hat{\mathbf l}=\frac{1}{\sqrt{p^2+q^2}}(p,q)$ is the unit vector.
In addition, by considering 
\begin{equation}
\begin{array}{ll}
& \gamma = \displaystyle\frac{ 4 b(\beta q-p)-2\alpha}{5 b(\beta^2+1)}, \\
& \nu= \displaystyle\frac{3 b (\beta^2+1)}{2 m_1 n_1},
\end{array}
\end{equation}
the solution (\ref{solA}) can be rewritten as
\begin{equation}\label{solAr}
\displaystyle u= 
\nu
\exp\big(-|{\mathbf l}| \hat{\mathbf l} \cdot {\mathbf x}\big)
\exp\big(\gamma({\mathbf k}\cdot{\mathbf x}-\alpha t)\big) 
\mathcal{P}\left(\frac{1}{\gamma}\exp\left({\frac{\gamma}{2}({\mathbf k}\cdot{\mathbf x}-\alpha t)}\right)+c_1, 0, c_2\right).
\end{equation}

\begin{figure}[h]
\centerline{\includegraphics[scale=0.4]{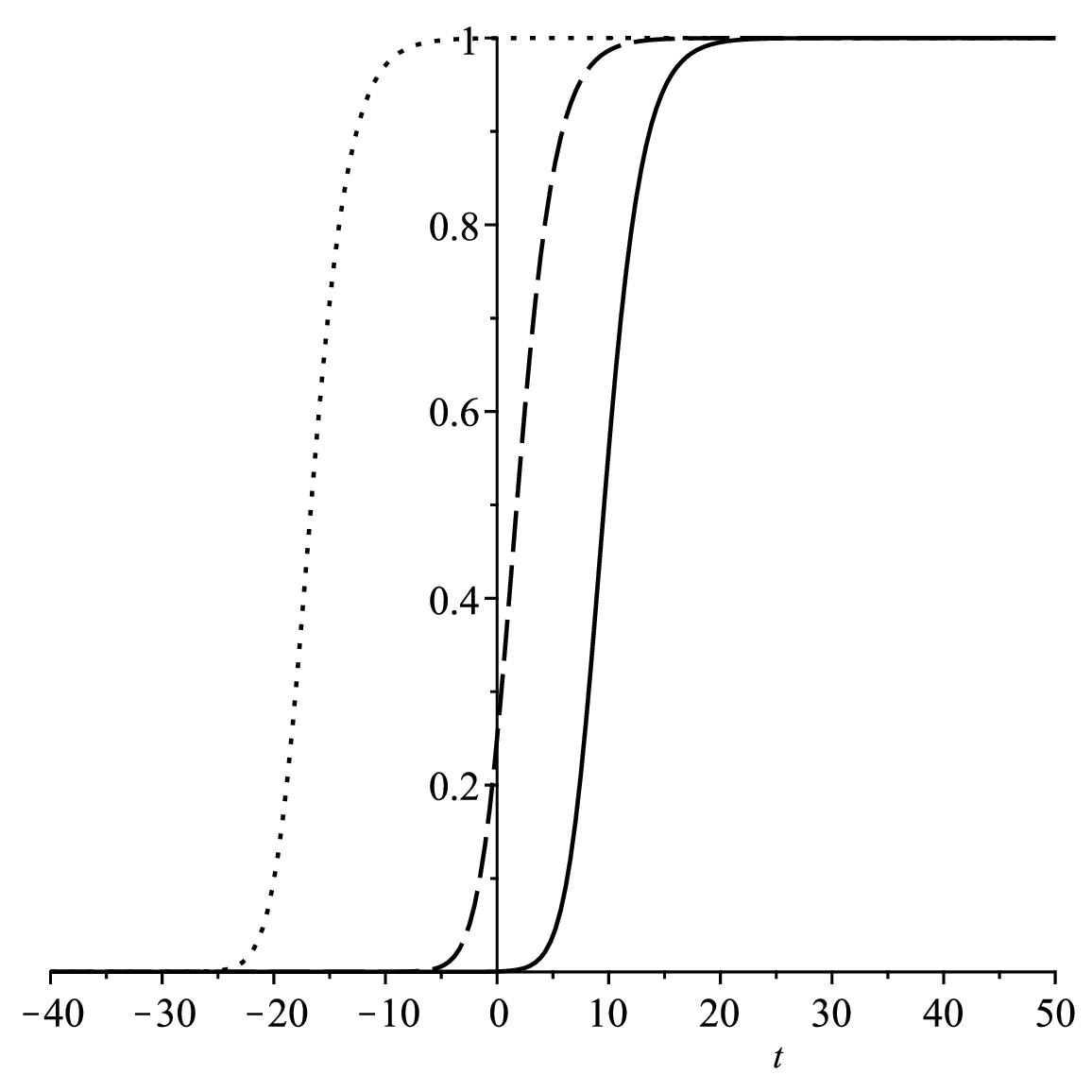}}
\caption{Shock solution (\ref{solAr}) at $(t,0,0)$ with $c_1=50$ (solid), $c_1=1$ (dash), $c_1=10^{-4}$ (dot). \label{fig1a}}
\end{figure}

\begin{figure}[h]
\centerline{\includegraphics[scale=0.4]{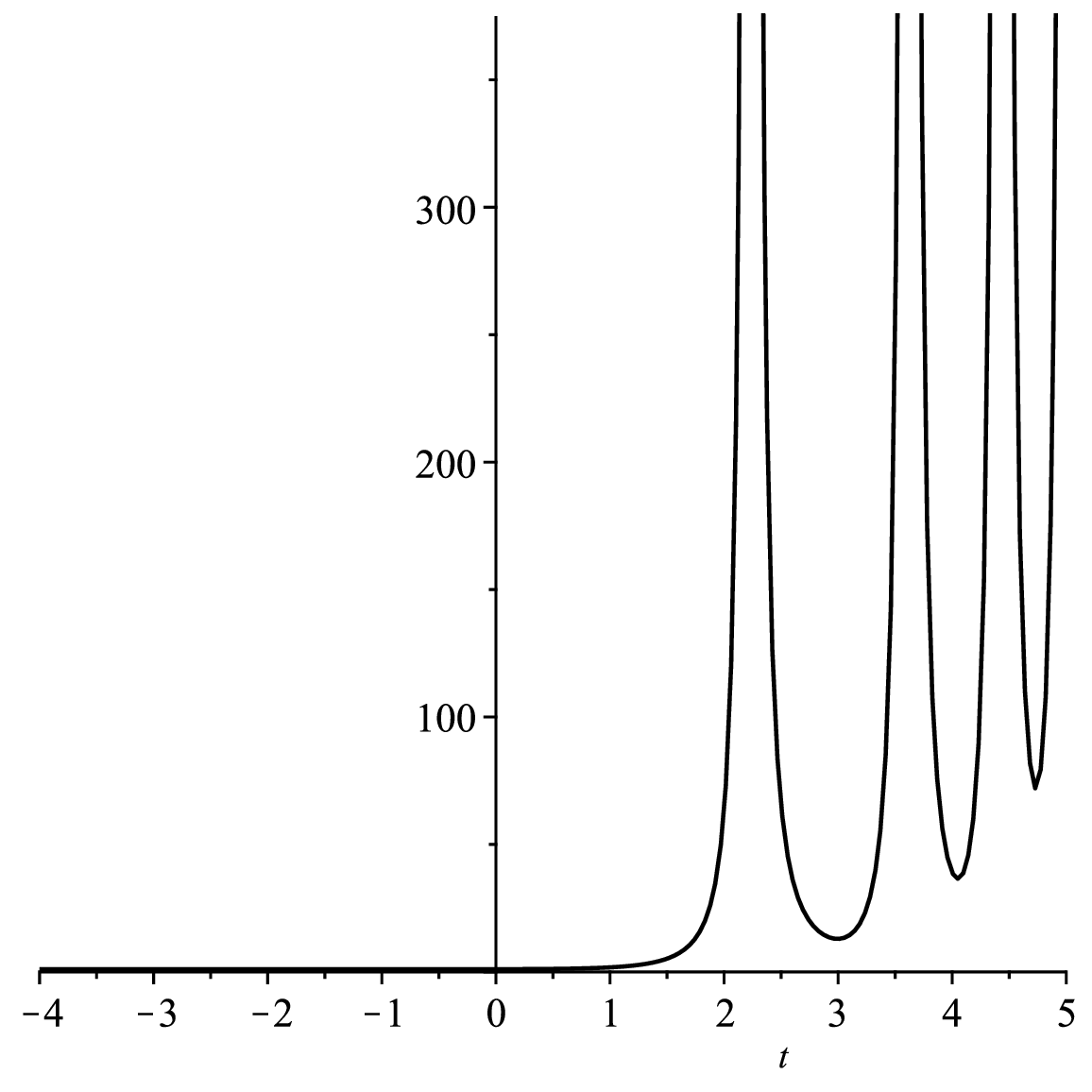}}
\caption{Solution (\ref{solAr}) at $(t,0,0)$ which blows-up in finite times. \label{fig2a}}
\end{figure}

\noindent
Because of the nature of the equation (\ref{edp}), the diffusion coefficient $b>0$, the growth rate $a$ could be positive or negative, and the parameters $\sigma$ and $\nu$ should be positives. When we consider $c_2=0$ and $c_1\neq0$ such that $\frac{1}{\gamma}\exp\big({\frac{\gamma}{2}({\mathbf k}\cdot{\mathbf x}-\alpha t)}\big)+c_1\neq0$, then solution (\ref{solAr}) is a shock transition when ${\mathbf x}$ is orthogonal to ${\mathbf l}$. 
For the values of the parameters $\gamma=1$, $\nu=1$, $\alpha=-1$, we show in Figure \ref{fig1a} the shock transition solution (\ref{solAr}) at $(t,0,0)$ for different values of $c_1$.
In other case, solution (\ref{solAr}) is the product of an exponential function and a real periodic function with period $c_2$, which blows-up in finite times. See Figure \ref{fig2a}.

\subsection{Reduction under the two-dimensional abelian subalgebra $\left\{ X_2, X_{2'} \right\}$ }

Consider $M(x)=m_1 e^{p x}$ and $N(y)=n_1 e^{q y}$ and the two-dimensional abelian subalgebra
\begin{equation}\label{algebra2}
\mathscr{B}:  \mbox{span} \left(  X_2, X_{2'}\right).
\end{equation}

By taking into account the two-dimensional abelian subalgebra (\ref{algebra2}), one can reduce (\ref{edp}) into a nonlinear first-order ODE. Let $X=X_2$ and $Y=X_{2'}$. From (\ref{conditions}), we obtain two independent invariants $z$ and $w$ given by
$$ \displaystyle z=t, \qquad w= e^{p x+q y} u.$$

\noindent Consequently, 
\begin{equation}\label{grinvsol2}
u= e^{-(p x-q y)} w \left(t \right),
\end{equation}

\noindent is a group-invariant solution. Taking into account the group-invariant solution (\ref{grinvsol2}), equation (\ref{edp}) reduces into the first-order ODE
\begin{equation}
\label{ode2} w'=-m_1 n_1 w^2+ (b(p^2+q^2)+a)w,
\end{equation}
whose general solution is given by
\begin{equation}
\label{sol2}  w(z)= \frac{b(p^2+q^2)+a}{c_1(b (p^2+q^2)+a) e^{-(b(p^2+q^2)+a)z}+m_1 n_1}.
\end{equation}

\noindent Undoing transformation (\ref{grinvsol2}), the general solution of the (2+1)-dimensional inhomogeneous Fisher-Kolmogorov equation (\ref{edp}) is
\begin{equation}\label{solB}
u(t,x,y)= \frac{ \left(b(p^2+q^2)+a\right) {e^{-px-qy}}}{c_1(b (p^2+q^2)+a) e^{-(b(p^2+q^2)+a)t}+m_1 n_1}.
\end{equation}
\\

\noindent
To interpret the previous solution, we rewrite it in a geometrical form. 
Since the coefficient of $u^2$ in equation (\ref{edp}) is given by $M(x)N(y)=m_1 n_1 e^{px+qy}$,
to interpret the solution,
it is natural to consider $px+qy= (p,q) \cdot (x,y) = {\mathbf l} \cdot {\mathbf x} = |{\mathbf l}| \hat{\mathbf l} \cdot {\mathbf x}$, where $\hat{\mathbf l}=\frac{1}{\sqrt{p^2+q^2}}(p,q)$ is the unit vector, and 
\begin{equation}
\begin{array}{ll}
& \gamma = b(p^2+q^2)+a  = b |{\mathbf l}|^2 + a, \\
& \nu= c_1\big(b (p^2+q^2)+a\big)  = c_1 \gamma,\\
& \sigma= m_1n_1.
\end{array}
\end{equation}
Hence, solution (\ref{solB}) can be rewritten as
\begin{equation}\label{solBr}
u= \frac{\gamma {e}^{\displaystyle -|{\mathbf l}|\ \hat{\mathbf l}\cdot {\mathbf x}}}{\displaystyle \nu {e^{-\gamma t}}+\sigma}.
\end{equation}
\\

\begin{figure}[h]
\centerline{\includegraphics[scale=0.4]{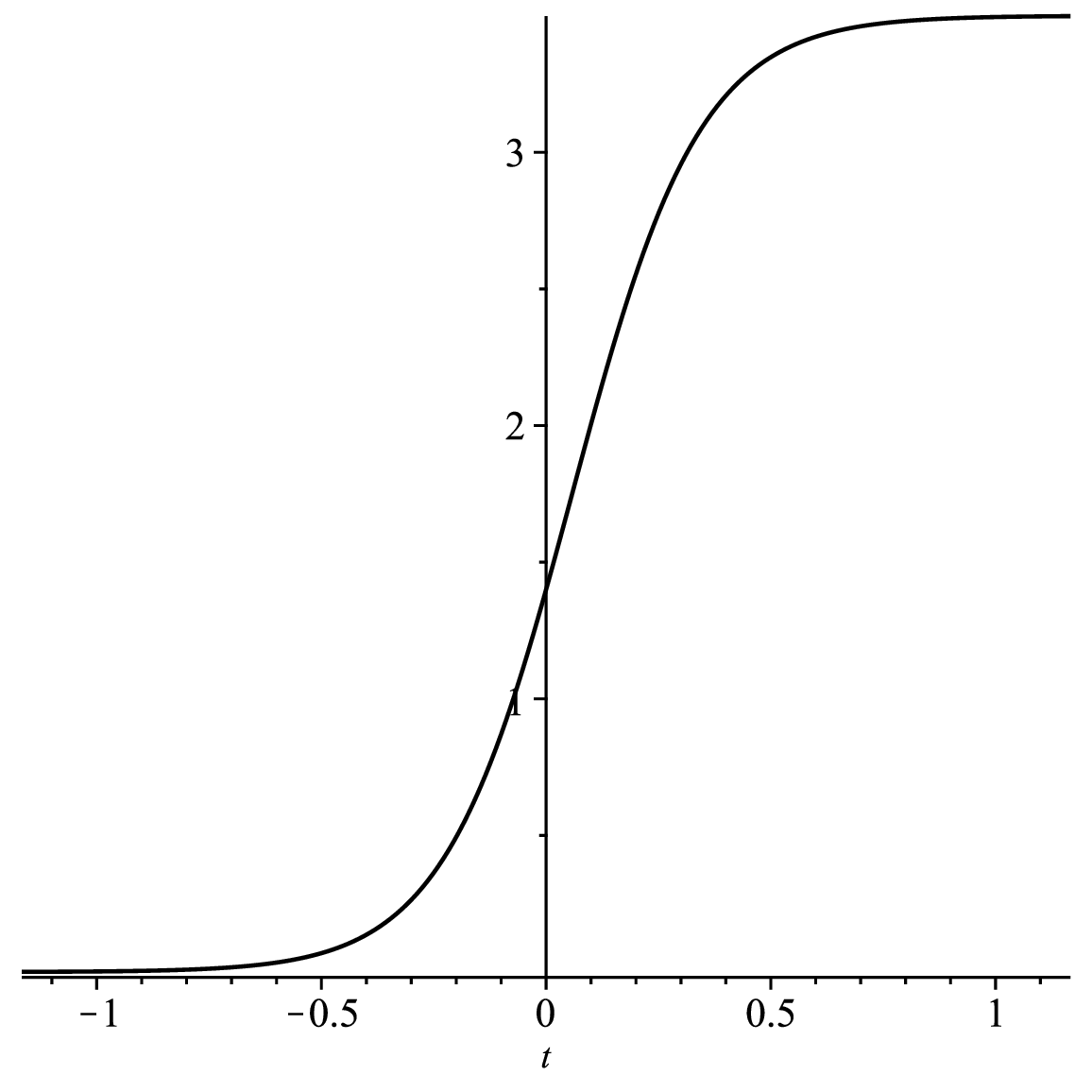}}
\caption{Shock transition solution (\ref{solBr}) $u(t,0,0)=\displaystyle \frac{7}{3 {e}^{-7t}+2}.$ \label{fig1b}}
\end{figure}

\begin{figure}[h]
\centerline{\includegraphics[scale=0.4]{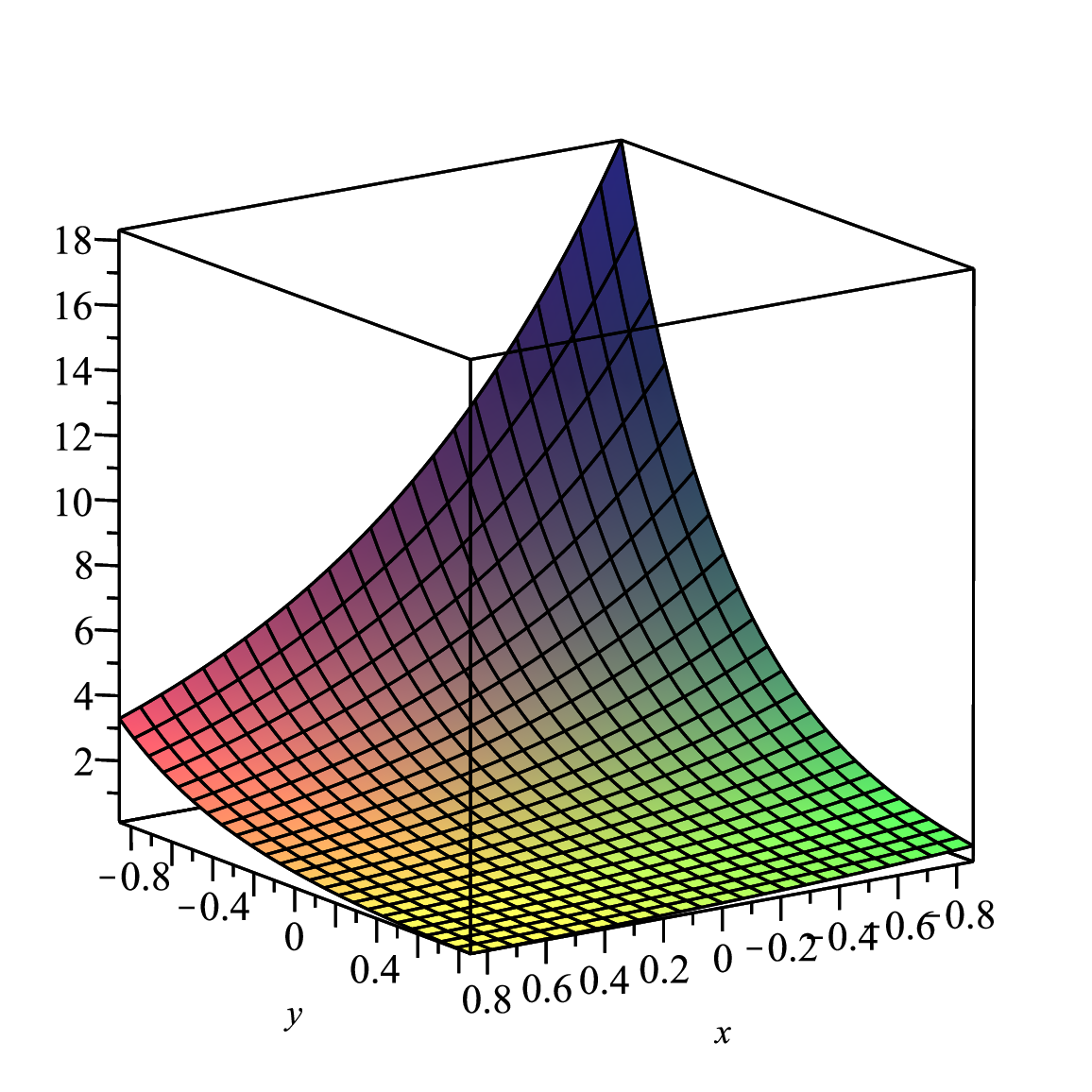}}
\caption{Spatial distribution of solution (\ref{solBr})) $u(0,x,y)=\displaystyle \frac{7}{5} {e}^{-x-2y}.$ \label{fig2b}}
\end{figure}

\noindent
Because of the nature of the equation (\ref{edp}), the diffusion coefficient $b>0$, the growth rate $a$ could be positive or negative, and also the arbitrary parameter $\nu$. We examine the case with $\nu>0$ and $a>0$. In Figures \ref{fig1b} and \ref{fig2b}, the solution (\ref{solBr}) with $a=2$, $b=1$, $p=1$, $q=2$, $\sigma=2$, $\nu=3$, is shown. Figure \ref{fig1b} shows the shock transition solution  (\ref{solBr}) at $(t,0,0)$, while Figure \ref{fig2b} shows the spatial distribution of solution (\ref{solBr}), which is an exponential sheet.

\section{Conclusions}
In this paper we have obtained all the Lie point symmetries admitted by the generalized Fisher-Kolmogorov equation in (2+1)-dimensions (\ref{edp}) with two exponential functions. By using the commutation of those point symmetries, two abelian subalgebras of dimension two have been considered. Hence, the (2+1)-dimensional inhomogeneous Fisher-Kolmogorov equation is reduced to ordinary differential equations. The new solutions obtained for (\ref{edp}) represent a shock transition, blows-up in finite times and exponential sheet behaviour. For future work, we plan to study the Lie symmetries of equation (\ref{edp}) for the case when $M(x)$ and $N(y)$ are arbitrary functions, and to use systematically all of the admitted Lie symmetries in order to get exact solutions.

\section*{ACKNOWLEDGMENTS}
The authors thank the reviewers for their careful reading of this article and their valuable remarks.
The authors express their sincerest gratitude to the Plan Propio de Investigaci\'on de la Universidad de C\'adiz MV2018-466.


\section*{CONFLICT OF INTEREST}
The authors declare no conflict of interest.

\bibliography{Bib_Icnaam2018v11}%

\end{document}